\documentclass[12pt,emtex]{article}
\usepackage{times}
\usepackage{graphicx}
\usepackage{amsfonts} %
\usepackage{amsmath} 
\usepackage{amssymb} %
\usepackage{amsthm}
\usepackage{epsfig}
\usepackage{latexsym}
\usepackage[matrix,curve]{xypic}
\usepackage[german]{babel} 
\newcommand{\C}{\Bbb{C}}

\newcommand{\Z}{\mathbf{Z}}

\begin{document}

\bigskip\noindent
\centerline{\bf Semisimple algebraic tensor categories}

\bigskip\noindent
\centerline{(Rainer Weissauer)}

\bigskip\noindent

\bigskip\noindent
For a field  $k$ a monoidal $k$-linear category abelian category $T$
is an abelian $k$-linear category with biadditive tensor functor $\otimes: T\times T \to T$,
$k$-linear and exact in each variable,  with associativity and commutativity constraints
and unit element $1_T$ satisfying the axioms ACU of [SR]. 
Then $T$ is called rigid, if every object $X$
has a dual $X^*$ with  morphisms $$\delta_X: 1_T \to X \otimes X^* \quad , \quad ev_X: X^* \otimes X \to 1_T$$
so that $(id_X \otimes ev_X)\circ ( \delta_X \otimes id_X) =id_X$ and 
$ ( \delta_X \otimes id_{X^*})\circ (id_{X^*} \otimes ev_X) =id_{X^*}$. See [CP] or [SR] for a detailed exposition.

\bigskip
Under the assumptions above, if $T$ is a small category such that
$End_T(1_T)\cong k$, the category $T$ is called a \lq categorie $k$-tensorielle \rq in [D].
If in addition $T$ is generated by one of its objects $V$ as a tensor category, such that
for some integer $N$ the lenght $l_T(V^{\otimes r})$ in $T$ is bounded by $N^r$ for all
$r$, the category $T$ will be called an algebraic tensor category over $k$. 

\bigskip
The typical example for an algebraic tensor category over $k$ (see [D], p.228)
is the category of finite dimensional $k$-linear $\varepsilon$-super representations $$ T = Rep_k({\bf G},\varepsilon)$$ of a super-affine
groupscheme ${\bf G}$ over $k$. The main result on algebraic tensor categories is
the following

\bigskip\noindent
{\bf Theorem 1}. ([D]) {\it Suppose $k$ is algebraically closed of characteristic zero.
Then any algebraic tensor category over $k$ is  of the form
$Rep_k({\bf G},\varepsilon)$.}

\bigskip
So let $k$ be algebraically closed of characteristic zero.
Under this assumption it is then interesting to know the cases where the category
$Rep_k({\bf G},\varepsilon)$ is a semisimple abelian category.
It is very easy to see that this only depends on the super-affine
groupscheme ${\bf G}$ and not on the additional twist $\varepsilon$.
In other words $Rep_k({\bf G},\varepsilon)$ is semisimple if and only if
the category $Rep_k({\bf G})$ of all $k$-linear finite dimensional super representations
of ${\bf G}$ is semisimple. More or less by definition $Rep_k({\bf G})$
coincides with the tensor category $CoRep_k(A)$ of $k$-finite dimensional
$A$-comodules, where $A$ is the super-affine Hopf algebra over $k$
defined the coordinate ring ${\cal O}({\bf G})$ of ${\bf G}$. If these categories are semisimple, we say ${\bf G}$ is reductive.

\bigskip
For a super-affine groupscheme ${\bf G}$ over a field $k$ of characteristic zero $k$ the reduced groupscheme of ${\bf G}$ is an algebraic group $G$ over
$k$. The left-invariant super derivations of the underlying
Hopf algebra $A$ corresponding to ${\bf G}$ define a finite dimensional Lie superalgebra $g=Lie({\bf G})$ over $k$. A Lie superalgebra $g$ over $k$ will be called reductive
if modulo its supercenter it is isomorphic to a direct sum of simple Lie superalgebras over $k$
of the classical types $A_n \ (n\geq 1), B_n \ (n\geq 3), C_n \ (n\geq 2), D_n\ (n\geq 3), E_6, E_7, E_8, G_2, F_4$ and of the orthosymplectic simple supertypes $BC_r \ (r\geq 1)$.
We then show 

\bigskip\noindent
{\bf Theorem 2}. {\it  ${\bf G}$ is reductive if and only its reduced group $G$ is a reductive
algebraic group over $k$ and its Lie superalgebra $Lie(\bf G)$ is  reductive over $k$.}

\bigskip
In particular ${\bf G}$ is reductive if and only if its connected component
${\bf G}^0$ with respect to the Zariski topology is reductive. In the connected
case we show that ${\bf G}$ is reductive if and only if etale unramified coverings are connected.

\bigskip
For the proof of theorem 2 we pass from super-affine groupschemes ${\bf G}$ over $k$ defined by their  super-affine Hopf coordinate algebra $A$ over $k$, to their associated supergroups $(G,g_-,Q)$. Here $G$ is the reduced group of ${\bf G}$.
The even part $g_+$ of $Lie({\bf G})=g_+\oplus g_-$ is the Lie algebra of $G$. The odd part $g_-$ is an algebraic $G$-module, and the Lie superbracket defines a $G$-equivariant symmetric map $Q: g_-\times g_- \to g_+$.  Together these data give rise to a triple $(G,g_-,Q)$ called a supergroup or a Harish-Chandra triple. For 
a suitable notion of representations for supergroups then the following holds

\bigskip\noindent
{\bf Theorem 3}. {\it The categories of $k$-finite dimensional super representations $Rep_k({\bf G})$ and $Rep_k(G,g_-,Q)$ are equivalent as algebraic tensor categories over $k$.} 

\bigskip
Theorem 3 allows us to reduce the proof of theorem 2 to the classical results
on the reductivity of semisimple Lie superalgebras obtained by Djokovic and Hochschild [DH].

\bigskip\noindent

\goodbreak
\bigskip\noindent
\centerline {\bf  Affine super Hopf algebras}

\bigskip\noindent
Let $k$ be field of $char(k)\neq 2$ and $A$ be a Hopf algebra with comultiplication, counit and antipode $(m_A^*,e^*_A,i^*_A)$ over the field $k$. Suppose $A$ is super-affine, i.e. suppose that as a ring A is a finitely generated super-commutative $k$-algebra such that $(m_A^*,e_A^*,i_A^*)$ are morphisms in the category $(salg)$ of super-commutative $k$-algebras.

\bigskip\noindent
{\it Remark}. The tensor product $\otimes^\varepsilon$ of the category $(salg)$ is the ordinary tensor product $\otimes_k$ except that it carries an induced grading with additional sign rules for certain structures like the tensor product of super $k$-algebras etc. For a detailed exposition of this we refer to [DM]. 

\bigskip
For the $\Z/2\Z$-grading $A=A_+\oplus A_-$ defined by the super structure the super-commutativity rule  $xy = (-1)^{\vert x\vert \vert y\vert} yx$ implies $x^2 =0$ for $x\in A_-$. 
Thus $A_-$ and the ideal $J$ generated by $A_-$ in $A$ are nilpotent. We call $J$ the super radical
of $A$. $J$ is a {\it Hopf ideal}, i.e. $i_A^*(J) \subset J$, $e_A^*(J)=0$ and   $$m_A^*(J) \subset J \otimes^\varepsilon A + 
A \otimes^\varepsilon J \ $$ as an immediate consequence of $$J= A_- + (A_-)^2\ $$ and $m_A^*(A_-)
\subset (A\otimes^\varepsilon A)_- \subset A_-\otimes^\varepsilon A + A\otimes^\varepsilon A_-$.
Surjective  Hopf algebra homomorphisms $\pi: A\to A'$ are in 1-1 correspondence with
Hopf ideals $I=Kern(\pi)$ of $A$. Since $A/J$ is even, the quotient $$\pi:A\to B=A/J$$ defines an commutative affine Hopf algebra quotient $B$ for which therefore $$G = Spec(B)$$ is a group scheme of finite type over $k$. We say $A$ is connected, if $G$ is connected in the Zariski topology.
Similar for the notion of being simply connected.
If $char(k)=0$, then $G$ is automatically reduced by a result of Cartier. In this case the super radical  $J$ is the nilradical of $A$.

\bigskip\noindent

\goodbreak
\bigskip\noindent
\centerline{\bf $A$-comodules}

\bigskip\noindent
An $A$-comodule $(V,\Delta_V)$ is a $k$-super vector space $V$ together
with a $k$-superlinear map
$$ \Delta_V: \ V \to V \otimes^\varepsilon A $$
satisfying the axioms (Modass) and (Modun) as in [S], p.30, i.e. the commutativity of
$$ \xymatrix{  V \ar[r]^-{\Delta_V}\ar[d]_{\Delta_V} & V\otimes^\varepsilon A \ar[d]^-{\Delta_V \otimes^\varepsilon id_A}
\cr  V\otimes^\varepsilon A \ar[r]^-{id_V\otimes^\varepsilon m_A^*}   &    V\otimes^\varepsilon A \otimes^\varepsilon A \cr } \quad  \quad 
\xymatrix{  V \ar[r]^-{\Delta_V}\ar[d]_{id_V} & V\otimes^\varepsilon A \ar[d]^{id_V \otimes^\varepsilon e_A^*}
\cr  V \ar@{=}[r]   &  V\otimes^\varepsilon k \cr }
$$
The notion of $A$-comodule only depends on the cogebra structure of $A$.
With the obvious notion of $A$-comodule homomorphism (see [S], p.31) the category of $A$-comodules is an abelian category. Any $A$-comodule is a union of its $k$-finite dimensional $A$-submodules.The category $CoRep_k(A)$ of $k$-finite dimensional $A$-comodules is a $k$-linear rigid abelian (monoidal) tensor category (see [CP], p.141). 

\bigskip\noindent
{\it Example a)}. $(A,m_A^*)$ itself is an $A$-comodule by the Hopf algebra axioms
$$ \xymatrix{  A \ar[r]^-{m_A^*}\ar[d]_{m_A^*} &  A\otimes^\varepsilon A \ar[d]^{m_A^* \otimes\varepsilon id_A}
\cr  A\otimes^\varepsilon A \ar[r]^-{id_A\otimes^\varepsilon m_A^*}  &  (A\otimes^\varepsilon A) \otimes^\varepsilon A \cr } \quad  \quad 
\xymatrix{  A \ar[r]^-{m_A^*}\ar[d]_{id_A} & A\otimes^\varepsilon A \ar[d]^{id_A \otimes^\varepsilon e_A^*}
\cr  A \ar@{=}[r]   &  A\otimes^\varepsilon k \cr }
$$

\bigskip\noindent
{\it Example b)}. For a $k$-super subvectorspace $V\subset A$ such that $ m_A^*(V) \subset V\otimes^\varepsilon A$ the restriction $\Delta_V = m_A^*\vert V$ defines an $A$-comodule $(V,\Delta_V)$, a subcomodule of $A$.

\bigskip\noindent
{\it Example c)}. Any Hopf algebra quotient $\pi: A\to B$ map makes $A$-comodules $(V,\Delta_V)$ into $B$-comodules $(V,\Delta)$ with respect to $$\Delta=(id_V\otimes^\varepsilon \pi)\circ \Delta_V : \quad V \to V\otimes^\varepsilon B \ .$$  
This is a consequence of $(\pi\otimes^\varepsilon \pi)\circ m_A^* = m_B^*\circ \pi$
and $e_A^*\circ \pi =e_B^*$. 

\bigskip\noindent

\goodbreak
\bigskip\noindent
\centerline{\bf Representations}

\bigskip\noindent
Suppose for an $A$-comodule $(V,\Delta_V)$ that the super vectorspace $V=k^{r\vert s}$ is finite dimensional with basis $e_i$ for $i=1...r+s$. Then $\Delta_V(e_i) = \sum_j e_j \otimes^\varepsilon f_{ji} $ for certain $f_{ji}\in A $. 
The axiom (Modass) implies $m_A^*(f_{ki}) = \sum_k f_{kj} \otimes^\varepsilon f_{ji}$. Thus the coefficients $f_{ji}$ define a homomorphism of super Hopf algebras
$$ {\cal O}(Gl(V)) \to A $$
from the super Hopf algebra $A'={\cal O}(Gl(V))$ of the general linear group of the super vector space $V$
to $A$.  Indeed as $k$-algebra $A'=k[X_{ij}, det_1^{-1},det_2^{-1}]$ is generated by elements $X_{kj}$ and the inverse of the determinants $det_1,det_2$ of the $X_{ij}$ for $i,j \leq r$ resp. $i,j > r$  
subject to the rule $m_{A'}^*(X_{ki})= \sum_k X_{kj} \otimes^\varepsilon X_{ji}$.
The elements $X_{ij}$ are even iff $i,j \leq r$ or $i,j > r$.
In other words,  this defines a super representation of $Spec^\varepsilon(A)$, i.e. a homomorphism of super group schemes
$$ Spec^\varepsilon(A) \to Gl(V) \ .$$
Conversely, it is easy to see that this defines a 1-1 correspondence between $k$-finite dimensional  $A$-comodules $V$
and finite $k$-linear dimensional super representations $V$ of the Lie super group scheme $Spec^\varepsilon(A)$. The category $Rep_k(A)$ of such  $k$-finite dimensional super representations of $Spec^{\varepsilon}(A)$ is an algebraic tensor 
category over $k$. The following is well known (see [D])

\goodbreak
\bigskip\noindent
{\bf Lemma 1}. {\it  
This correspondence induces a tensor-equivalence between the algebraic tensor categories  $CoRep_k(A)$ and $Rep_k(A)$ over $k$.}

\bigskip\noindent

\goodbreak
 \bigskip\noindent
 \centerline{\bf The functor of invariants $V \mapsto V^G$}
 
 \bigskip\noindent
For the $k$-groupscheme $G=Spec(B)$ consider the left-exact functor $$ V \ \ \mapsto \ \ V^{G}= Hom_{B-comod}(k,V) \ $$
from the category of $B$-comodules to the category of $k$-vectorspaces.
The $k$-vectorspace  $V^G \subseteq V$ can be identified with the maximal trivial $B$-subcomodule of $V$ of all elements $v$ in $V$  for which
 $$\Delta_M(v) = v\otimes 1_B\ .$$ 

\bigskip\noindent
We say a $B$-comodule $V$ is free, if it isomorphic to a $B$-comodule of the form $V=V_0\otimes B = B^d$. Here $V_0$ is a $k$-vectorspace and $d=dim_k(V_0)$. 
$B$-comodules will be called almost
free, if they have a finite filtration by $B$-subcomodules whose sucessive quotients are free
$B$-comodules.
Notice
$B^G = k\cdot 1_B$, since  $v= ( e_B^*\otimes id_B )(m_B^*(v)) =
 ( e_B^*\otimes id_B )(\Delta_V(v)) = (e_B^*(v)\otimes 1_B) \in k\cdot 1_B $
 for $v \in B^G$. Hence for free $V=V_0 \otimes B$
$$   (V_0 \otimes B)^G = V_0 \ .$$

\bigskip\noindent 
Using bar-resolutions (see [DG], p.233ff) one can define derived  functors $H^i(G,-)$
such that $H^0(G,V) = V^G$.
In other words
a short exact sequence of $B$-comodules gives rise
to a long exact sequence of $k$-vectorspaces using the derived  functors $H^i(G,-)$. 
By [DG], lemma 3.4 $$H^i(G,B)=0 \quad , \quad i\geq 1$$ for any free $B$-comodule.  Obviously $H^1(G,V)=0$ for almost free $B$-comodules $V$. 
Hence 

\bigskip\noindent
{\bf Lemma  2}. {\it On the Grothendieck group of almost free $B$-comodules $V$}
$$  rang_k(V) = dim_k(V^G) $$
{\it defines a homomorphism .}

\bigskip\noindent

\bigskip\noindent

\goodbreak
\bigskip\noindent
\centerline{\bf The Hopf ideals defined by $J$}

\bigskip\noindent
Let $A$ be a super-affine Hopf algebra over $k$.
Then its super radical  $J$ is generated as an $A$-module by finitely many elements in $A_-$.
If $J$ is generated by $s$
elements then it is easy to see that $J^{s+1}=0$. 
Hence there exists a finite descending filtration by $A$-right (and left) ideals
$$ 0 \subset J^s \subset J^{s-1} \subset .. \subset J^2 \subset J \subset A $$
whose sucessive quotients $$ V_i=J^i/J^{i+1}$$ are right (and left) $B=A/J$-modules. 
Although the $J^i$ are not $B$-modules  a priori, they are
$B$-subcomodules  of the $B$-comodule $(A,\Delta)$ with structure map $$\Delta= (id_A\otimes^\varepsilon \pi)\circ m_A^*$$ using the examples a), b) and c) above. There is a commutative diagram
$$ \xymatrix{ A  \ar[r]^-{m_A^*} &  A\otimes^\varepsilon A \ar[r]^-{id_A\otimes^\varepsilon \pi} &  A\otimes^\varepsilon B \cr
J^i \ar[rr]^{\Delta} \ar@{^{(}->}[u]  &   &  J^i \otimes^\varepsilon B \ar@{^{(}->}[u]\cr} $$
since the image of $ m_A^*(J^i) \subset m_A^*(J)^i \subset (J \otimes^\varepsilon A + 
A \otimes^\varepsilon J)^i \subset \sum_{a+b=i} J^a \otimes^\varepsilon J^b $ in $A \otimes^\varepsilon B$, under
$id_A\otimes^\varepsilon \pi$, is contained in 
$ J^i \otimes^\varepsilon B$. Thus $J^i$ becomes a $B$-comodule. The $V_i$ then are quotient $B$-comodules of the $J^i$
in the obvious way.

\bigskip\noindent
{\bf Lemma 3}. {\it $V_i  \cong B^{d_i}$ is finite free both as a $B$-right module and a $B$-comodule.}

\bigskip\noindent
{\it Proof}. The $k$-linear structure map $\Delta: A \to A\otimes^\varepsilon B$ of the $B$-comodule $A$ is  $A$-linear in the following sense: For $a\in A$ and $x\in A$
of course $x\cdot a \in A$. Since $\pi$ is $A$-linear $$ \Delta(x\cdot a) = (id_A\otimes^\varepsilon \pi)\bigl(m_A^*(x)\cdot m_A^*(a)\bigr) =  \Delta(x) \bullet
 m_A^*(a) $$
where $A\otimes^\varepsilon B$ is viewed as a $A\otimes^\varepsilon A$-right module
in the obvious way. In other words   
 $ m_A^*(a) = \sum a_\nu \otimes a'_\nu $ acts on $y= a\otimes^\varepsilon b$
 via $y \bullet  m_A^*(a) = \sum_\nu (-1)^{\vert a_\nu \vert \vert b\vert}a \cdot a_\nu \otimes^\varepsilon b \cdot a'_\nu 
 $.
Since $\Delta(J^i) \subset (J^i)$ the map $\Delta$ induces a quotient map
$$ \Delta_V: \  V \to V \otimes^\varepsilon B \ $$
on $V=V_i=J^i/J^{i+1}$ making it to a $B$-comodule. The right action of $A$ on $V$ factors over the quotient
ring $B$. Similarly the right action of $A\otimes^\varepsilon A$ on $V\otimes^\varepsilon B$
factors over the quotient ring $A/J \otimes^\varepsilon B = B \otimes^\varepsilon B$, so that now (*)
$$ \Delta_V(x \cdot b) = \Delta_V(x) \bullet m_B^*(b) \ $$ is obvious: The composition of $m_A^*$ with the projection
$A \otimes^\varepsilon A \to A/J \otimes^\varepsilon B$
is equal to $m_B^* \circ \pi$.

\bigskip
It is the property (*) which makes the right $B$-module and right $B$-comodule $V$ into a $B$-right {\it Hopf module} in the sense of [S], p.83. Since $B$ is an ordinary Hopf algebra
we can immediately apply [S], theorem 4.1.1. It states that
$$      M \cong M^G \otimes B = B^d \quad , \quad  d=dim_k(M^G) $$
as a Hopf right $B$-module and comodule for any Hopf right $B$-module and co-module $M$.
Applied for $M=V$ we now use the fact that $J$, hence also $V$, are finitely generated $B$-right modules. Hence $d=d_i < \infty$ in our case. This proves our claim. QED

\bigskip
Therefore $A$ is an almost free $B$-comodule. By lemma 2 this implies 

\bigskip\noindent
{\bf Corollary 1}. {\it   $\ dim_k(A^G) = \sum_{i=0}^s d_i\ $ for $d_i=rank_B(J^i/J^{i+1})$.}

\bigskip\noindent
{\it Remark}. We will see later
in corollary 5 that for $A$ affine super group scheme over $k$  we have $d_i= {s\choose i} $. This will imply $$ dim_k(A^G) = 2^s \ .$$

\bigskip\noindent
{\bf Lemma 4}. {\it $A^G$ is a finite dimensional $k$-subalgebra of $A$}.

\bigskip\noindent
{\it Proof}. $\Delta(v) = v\otimes 1_B$ and $\Delta(v') = v'\otimes 1_B$
imply $\Delta(v\cdot v') = \Delta(v)\cdot \Delta(v')= (v\otimes 1_B)\cdot (v'\otimes 1_B)
= (v\cdot v')\otimes 1_B$. QED

\bigskip\noindent

\goodbreak
\bigskip\noindent
\centerline{\bf Superderivations}. 

\bigskip\noindent
Let $(A,m^*_A,e^*_A,i^*_A)$ be an super-affine Hopf algebra over $k$.
Let $m=m(A)=kern(e_A^*)$ be the maximal ideal  of $A$ at the identity.
Tanget vectors $X\in (m/m^2)^*_{\pm}$ extend to even or odd
$k$-linear superderivations $d_X:A\to k$ by composing $X:  m/m^2 \to k$ with the projection $A=k\cdot 1\oplus m \to m \to m/m^2$. 
Define $k$-super derivations 
 $$ D_X:A\to A $$ by the commutative diagram
$$\xymatrix{ A \ar[r]^-{m_A^*}\ar[dr]_-{D_X} &  A \otimes^\varepsilon A \ar[d]^{id_A\otimes^\varepsilon d_X}\cr &   A \cr }$$
Then  $d_X=e_A^* \circ D_X $ by definition. The  $k$-superderivations $D_X:A\to A$
so constructed are left-invariant, i.e. for $D=D_X$ there exists a commutative diagram
$$\xymatrix{ A \ar[r]^-{m_A^*}\ar[d]_-{D} &  A \otimes^\varepsilon A \ar[d]^{id_A\otimes^\varepsilon D}\cr 
A \ar[r]^-{m_A^*}  &   A\otimes^\varepsilon A \cr }$$
by the coassociativity law
$(id_A \otimes^\varepsilon m_A^*) \circ m_A^* = (m_A^*\otimes^\varepsilon id_A) \circ m_A^*  $
and $A \otimes^\varepsilon (A \otimes^\varepsilon A) = (A \otimes^\varepsilon A) \otimes^\varepsilon A$. Indeed, if  we apply
$id_A \otimes^\varepsilon (id_A \otimes^\varepsilon d_X) = (id_A \otimes^\varepsilon id_A )\otimes^\varepsilon d_X$ on the left side
of the coassociativity law, this becomes $(id_A \otimes^\varepsilon D_X) \circ m_A^*$.
On the right side of the coassociativity law it becomes $   m_A^* \circ  D_X $.

\bigskip\noindent
{\bf Lemma 5}. {\it There exists a canonical isomorphism $X\mapsto D_X$ of $k$-vectorspaces
$$ (m/m^2)^*  \to Lie(A) $$
 between the tangent space at the identity element and the $k$-vector space $Lie(A)$ of all
 left-invariant $k$-superderivations of $A$.}

\bigskip\noindent
{\it Proof}. The inverse map is $Lie(A)\ni D\mapsto d=e_A^*\circ D$. Since 
$d:A \to k$ is a $k$-superderivation, it must vanish on $m^2$ and on $k\cdot 1$.
Hence $d=d_X$ for some $X\in (m/m^2)^*$. The left-invariant $k$-superderivation $D:A\to A$ is uniquely determined by
its restriction $d=e_A^* \circ D$, since $d$  determines $D$ via the right vertical arrow $id_A \otimes^\varepsilon d$ of the composed commutative diagrams
$$\xymatrix{ A \ar[r]^-{m_A^*}\ar[d]_-{D} &  A \otimes^\varepsilon A \ar[d]^{id_A\otimes^\varepsilon D} \cr 
A \ar[r]^-{m_A^*}  \ar[dr]_{id_A}&   A\otimes^\varepsilon A \ar[d]^{id_A\otimes^\varepsilon e^*} \cr &  A \cr}$$

\bigskip\noindent
For any Hopf ideal $J$ of $A$ with quotient map $A\to B=A/J$ the cotangent space $m(A)/m(A)^2$ surjects onto the cotangent space
$m(B)/m(B)^2$. Hence $Lie(B)$ injects into $Lie(A)$. QED

\bigskip\noindent 
{\bf Lemma 6}. {\it The image of the natural injection $Lie(B)\hookrightarrow Lie(A)$  is the space of left-invariant $k$-derivations $D$ of $A$ (as in the last lemma)
with the property $D(J)\subset J$. }

\bigskip\noindent
{\it Proof}. Such $D$ induce left-invariant derivations on the quotient $B=A/J$.
So it suffices that $X\in Lie(B)$
implies $D_X(J) \subset J$.
For this
 let $x:A\to B$
be the quotient map with kernel $J$, considered as a $B$-valued point of $A$. For $f\in A$ by definition  $D_X(f)(x) = (x\otimes^\varepsilon d_X)(m_A^*(f))$.  Now $(x\otimes^\varepsilon d_X)(m_A^*(f)) \subset (x\otimes^\varepsilon d_X)(A \otimes^\varepsilon J + J \otimes^\varepsilon A)$ for $f\in J$ since $J$ is a Hopf ideal.  But $x(J)=0$. On the other hand
$d_X(J)=0$ for $X\in Lie(B)$, since $Lie(B)$ is the space of linear forms $d_X:m(A)/m(A)^2 \to k$ trivial on the image of $J$. Hence $D_X(f)(x)=0$ or $D_X(f)\in J$. QED

\bigskip\noindent
{\it The Lie algebra}. The supercommutator $[D,D'] = D\circ D' - (-1)^{\vert D\vert \vert D'\vert}D'\circ D$ of two $k$-super derivations $D,D'$ is a $k$-superderivations. Since the super commutator
of left-invariant derivations is left invariant, the finite dimensional $k$ super-vector space
  $$g=Lie(A) = (m(A)/m(A)^2)^* $$ defined by $(A,m_A^*,e_A^*)$
 becomes a  Lie $k$-superalgebra with $g\cong k^{r\vert s}$ as a super vectorspace
 $$g=g_+\oplus g_- \ .$$
 
\bigskip\noindent
{\it The super radical $J$}. Notice $J=A_-+(A_-)^2 $ implies  $J^2=(A_-)^2 + (A_-)^3$. Hence  the quotient $J/J^2 = A_-/(A_-)^3$ is odd.
Since $J$ is nilpotent, we have $J\subset m(A)$ and
$m(B) = m(A)/J$. Clearly the quotient $A_-/(A_-)^3=J/J^2 \to J/(J \cap m(A)^2)$ again is odd. 
Since $B$ is even, also $m(B)/m(B)^2$ is even with $g_+=Lie(B)=(m(B)/m(B)^2)^*$ even.
Hence the exact sequence
$$ 0 \to J/(m(A)^2\cap J) \to m(A)/m(A)^2 \to m(B)/m(B)^2 \to 0 $$
gives rise to a splitting of the super-vectorspace $Lie(A)$ with $Lie(G)$ even
$$  Lie(A) = Lie(G) \oplus g_- \ $$
and with $g_- \cong (J/(J \cap m(A)^2))^*$ odd.

\bigskip\noindent
Fix a basis $\tilde \theta_i$ of $(V_1)^G = (J/J^2)^G$ and representatives 
$\theta_i\in J^G$ of the elements $\tilde \theta_i$.  Then
$J/J^2 = \bigoplus_{i=1}^{d_1} \tilde\theta_i \cdot B$ as a $B$ right-module.
Consider the exact sequence of odd $k$-vectorspaces
$$ 0 \to K \to J/J^2 \to J/(J\cap m(A)^2) \to 0 \ .$$
We claim
$  K =  \bigoplus_{i=1}^{d_1} \tilde\theta_i \cdot m(B) $.

\bigskip
Since $\theta_i \in J \subset m(A)$, the right hand side is contained in $K$.
Conversely elements $k\in K$ have odd representatives $x$ in $J\cap m(A)^2$, or hence
in $A_- \cap m(A)^2$. Notice 
$A_- \cap m(A)^2 = (A_- \cap m(A))  (A_+ \cap m(A))$
by a case by case verification and the definition of the super graded ring structure on $A$. 
Since $m(A)\cap A_- \subset J$, hence $A_- \cap m(A)^2 
 \subset J\cdot (A_+ \cap m(A)) $. As $m(A)$ acts on $J/J^2$ via its quotient $m(B)$
 therefore the image $k$ of $x$ is contained in $\bigoplus_{i=1}^{d_1} \tilde\theta_i \cdot m(B)$.  This proves the claim. As a consequence 
$$  (J/J^2)^G \ = \ \bigoplus_{i=1}^{d_1}\  \tilde\theta_i \cdot k \ \cong  \  (J/J^2)\big/K \ \cong \ J/(J\cap m(A)^2) \ \cong \ (g_-)^* \ .$$
Together with the lemma 6 this implies

\bigskip\noindent
{\bf Corollary 2}. {\it The left-invariant derivations $D_X$ for $ X\in Lie(G) \subset Lie(A)$ respect the  exact sequence defined by the super radical $J$
$$ \xymatrix{  0 \ar[r] & J \ar[r] & A \ar[r] &  B \ar[r] & 0 \cr} \ .$$
A left-invariant superderivation $D_X\in Lie(A)$ preserves the super radical $J$ if and only if
$X \in Lie(G)$. Furthermore $$ dim_k(g_-)=rank_B(J/J^2)=d_1 \ .$$}

\bigskip\noindent
{\it Homomorphisms}. A homomorphism $\Phi^*: A \to A'$ between super-affine $k$-Hopf algebras  induces a map between the tangent spaces at the identity element, hence a $k$-linear map
$$ Lie(\Phi): \ Lie(A') \to Lie(A) \ .$$ 
$Lie(\Phi)$ is a homomorphism of $k$-super Lie algebras, since
$\Phi^* \circ D_X = D_{X'} \circ \Phi^*$
for $X'=Lie(\Phi)(X)$. [Reduce to $(\Phi^*\otimes^\varepsilon\Phi^*) \circ (id \otimes^\varepsilon d_X) = 
(id \otimes^\varepsilon d_{X'})\circ \Phi^*$, hence to $\Phi^* \circ d_X = d_{X'}$.]

\bigskip\noindent
{\it Adjoint action}. The interior automorphism $\Phi^*=(Int_x)^*$ defined by a $k$-valued point
of $Spec(A/I)$ induces a Lie algebra homomorphism $Ad(x)=Lie(Int_x)$
from $Lie(A)$ to $Lie(A)$.
Obviously $Ad(x)\circ Ad(y) = Ad(xy)$. Hence $Ad(x)$ defines a $k$-linear 
representation on $Lie(A)$ of the underlying algebraic group $G$ 
$$  Ad:\ G(k) \to   Gl_k(Lie(A)) \ .$$ 
This adjoint action respects the super structure, hence decomposes into 
representations $Ad_\pm$ of $G$ on $g_+$ and $g_-$ respectively. $Ad_+$ is the usual adjoint action of $G(k)$ on its Lie algebra $g_+=Lie(G)$. 

\bigskip\noindent

\goodbreak
\bigskip\noindent
\centerline{\bf Left versus right} 

\bigskip\noindent
Similar to left-invariant superderivations 
define right-invariant superderivations of a Hopfalgebra $A$. The Lie superalgebra of the left-invariant and right-invariant superderivations are isomorphic
(use the antipode).
Left-invariant superderivations $D$ and right-invariant superderivations $D'$ of $A$ supercommute. Use $$(-1)^{\vert D\vert \vert D'\vert}m_A^*(DD'x)=(D'\otimes^\varepsilon D)(m_A^*(x)) =m_A^*(D'Dx)\ $$ to show
that their supercommutator $[D,D']$ is a derivation with $m_A^*([D,D'](x))=0$.
Hence $[D,D']=0$ by applying the counit $e_A^*$.

\bigskip\noindent
{\bf Lemma 7}. {\it For quotients $B=A/I$ by a Hopf ideal $I$ and $X\in Lie(B) \subset Lie(A)$ the left-invariant superderivations $D_X$ of $A$  preserve $B$-subcomodules $V$ of $A$}.


\bigskip\noindent
{\it Proof}. The commutative diagram
$$\xymatrix{ A \ar[r]^-{m_A^*}\ar[d]_-{D_X} &  A \otimes^\varepsilon A \ar[d]^{id_A\otimes^\varepsilon D_X} \ar[r]^-{id_A\otimes^\varepsilon\pi} &  A \otimes^\varepsilon B \ar[d]^{id_A\otimes^\varepsilon D_X}\cr 
V \ar[r]^-{m_A^*}  &   A\otimes^\varepsilon A \ar[r]^-{id_A \otimes^\varepsilon \pi}  &   A\otimes^\varepsilon B \cr }$$ 
for $D_X$ and $X\in Lie(B) \subset Lie(A)$ gives $\Delta\circ D_X = (id\otimes^\varepsilon D_X)\circ \Delta$ for the structure map $\Delta= (id_A\otimes^\varepsilon \pi)\circ m_A^*$ of the $B$-comodule $A$. 

\bigskip\noindent
Next notice $(id_A\otimes e_B^*)\circ \Delta = id_A $ and the commutative diagram
$$\xymatrix{ A \ar[r]^-{m_A^*}\ar@{=}[d] &  A \otimes^\varepsilon A \ar[d]^{id_A\otimes^\varepsilon e_A^*} \ar[r]^-{id_A\otimes^\varepsilon\pi} &  A \otimes^\varepsilon B \ar[d]^{id_A\otimes^\varepsilon e_B^*}\cr 
A \ar@{=}[r] &   A \ar@{=}[r]  &   A \cr }
$$
 Hence
for $v\in V\subset A$ and $\Delta(v)=\sum_i v_i \otimes^\varepsilon b_i$ with $v_i\in V, b_i\in B $
the element $$ D_X(v)=(id_A\otimes^\varepsilon e_B^*)\circ \Delta(D_X(v))$$ is $(id_A\otimes^\varepsilon e_B^*)(\sum_i v_i \otimes^\varepsilon D_X(b_i)) )= \sum_i v_i \cdot d_X(b_i) $
using left-equivariance of $D_X$ as in first diagram above. Thus $D_X(v) \in V$ and $D_X(V)\subset V$. QED. 

\bigskip
For $\Delta(v)=v\otimes^\varepsilon 1_B$ in particular  $D_X(v)=0$, since $d_X(1_B)=0$.

\bigskip\noindent
{\bf Corollary 3}. {\it  $D_X(A^G)=0$ for all $D_X, X\in Lie(G)$ where $G=Spec^\varepsilon(A/I)$. }

\bigskip\noindent
{\bf Corollary 4}. {\it  $A^G$ is stable under all right-invariant superderivations $D'_X$
in  $ Lie(A)$. }

\bigskip\noindent
{\it Proof}. This follows from the commutative diagram 
$$\xymatrix{ A \ar[r]^-{m_A^*}\ar[d]_-{D'_X} &  A \otimes^\varepsilon A \ar[d]^{D'_X\otimes^\varepsilon id_A} \ar[r]^-{id_A\otimes^\varepsilon\pi} &  A \otimes^\varepsilon B \ar[d]^{id_A\otimes^\varepsilon D'_X}\cr 
V \ar[r]^-{m_A^*}  &   A\otimes^\varepsilon A \ar[r]^-{id_A \otimes^\varepsilon \pi}  &   A\otimes^\varepsilon B \cr }$$ 
which implies $ \Delta \circ D'_X =  (D'_X \otimes^\varepsilon id_B) \circ \Delta  $ 
for the structure map $\Delta$ of the $B$-comodule $A$.
For $v \in A^G$ by definition $\Delta(v) = v\otimes^\varepsilon 1_B$. Hence
$\Delta(D'_X(v))= (D'_X \otimes^\varepsilon id_B) \circ \Delta(v)= D'_X(v)\otimes 1_B$. This shows $D'_X(A^G)\in A^G$. QED

 \bigskip\noindent
\centerline {\bf The subring $A^G \subset A$}  

\bigskip\noindent
For the super radical $J$ of $A$ put $G=spec(B)$ and $B=A/J$ as before.
 Since $(J/J^2)$ is odd and almost free, the quotient map $(J_-)^G \to (J/J^2)^G$ is surjective so that we can choose representatives $\theta_1,..,\theta_s \in (J_-)^G$ of a $k$-basis
in $(J/J^2)^G$ so that the $\theta_i$ are also a $B$-basis of $J/J^2$ by lemma 2. Then by recursion modulo the $J^n$
$$  J = \theta_1 \cdot A + \cdots + \theta_s \cdot A \ .$$

\bigskip\noindent 
The $\theta_i$ are odd. Hence  by supercommutativity $$\theta_i \theta_j = - \theta_j \theta_i\ .$$  For $I\subset \{1,..,s\} $ define $\theta_I = \theta_{i_1} \cdots \theta_{i_n}$ if $I=\{i_1,..,i_n\}$ and $i_1 < ... < i_n$.
With these notations $J^n$ is generated as an $A$-right module by the
 $\theta_I$ with $\vert I\vert = n $. Hence for the elements $\tilde\theta_I = \theta_I$ mod
$J^{n+1}$ in  $(J^n/J^{n+1})^G$ we get $$ J^n/J^{n+1} = \sum_{\vert I\vert =n} \tilde\theta_I \cdot B  \ .$$
We may replace by a $B$-right linear independent subset of $T_n$ of the set of all the $\tilde\theta_I$, since we already
know that $J^n/J^{n+1}$ is a free $B$-right module generated by a $k$-basis of $(J^n/J^{n+1})^G$. Therefore
$$ J^n/J^{n+1} = \bigoplus_{I \in T_n} \tilde\theta_I \cdot B  \ $$
and
$$ (J^n)^G /(J^{n+1})^G \ \cong \ (J^n/J^{n+1})^G \ \cong \ k^{\# T_n}   \ .$$
Since $\theta_I \in A^G$, recursively now any element in $A^G$ can be written as a
superpolynom in the elements $\theta_1,...,\theta_s$ by induction modulo the  $A^G\cap J^n =(J^n)^G$.This defines a surjective
$k$-algebra homomorphism $  f:  S^\varepsilon(k^{0\vert s}) \to A^G$ mapping the generators of the superpolynomial ring $S^\varepsilon(k^{0\vert s})$ to the $\theta_i$.

\bigskip\noindent
{\bf Lemma 8}. {\it $A^G$ is a superpolynomial ring $S^\varepsilon(k^{0\vert s})$ over $k$ in the odd variables $\theta_i$}.

\bigskip\noindent
{\it Proof}. Recall $(J/J^2) \cong (g_-)^*$.
This means that
we can find $s$ odd right-invariant superderivations $D'_i$ 
in $g_-\subset Lie(A)$ such that $e_A^*(D'_i(\theta_j))= d'_i(\theta_j)
= \delta_{ij}$ in $k$. In other words
$$  D'_i(\theta_j) \ \equiv \  \delta_{ij} \ \mod m(A) \ .$$
Since $D'_i(A^G) \subset A^G$ and 
since $  m(A) \cap A^G = J^G $  $$D'_i(\theta_j) = \delta_{ij} + Q_{ij}(\theta)  $$
for certain super polynomials $Q_{ij}$ in the variables
$\theta_i$, whose minimal nonvanishing Taylor coefficient has degree $\geq 1$.
Suppose $P\neq 0$ is an element in $I=Kern(f)$ with minimal nonvanishing Taylor coefficient  
say of degree $d$, such that this $d$ is minimal among all $0\neq P\in I$.
If $d=0$, then $P$ is a unit in the superpolynomial ring and the quotient $A^G$ would 
be zero in contradiction to  $1_A\in A^G$. Hence $d>0$.
Let $\theta_i$ be a variable which occurs nontrivially in the Taylor coefficient
of $P$ of degree $d$. Then apply the derivative $D'_i(P)$.
Obviously $D'_i(P)$ has a nonvanishing Taylor coefficient of degree $d-1$.
On the other hand $D'_i(I) \subset I$, hence $D'_i(P)\in I$. This gives a
contradiction unless the kernel vanishes $I=0$. QED

\bigskip
Then by an obvious counting argument lemma 8 implies

\bigskip\noindent
{\bf Corollary 5}. {\it $d_n = \# T_n = {s \choose n}$ for all $n$}.

\bigskip\noindent
{\it Choice of bases}. Up to a scalar  $\eta = \theta_I$ for $I=\{1,..,s\}$  is independent
of the choice of the basis $\theta_i$, since it is a generator of the one dimensional
$k$ vectorspace $(J^s)^G$. Hence $\eta$  
is an eigenvalue of the right-invariant operators $D'\in Lie(G)$ corresponding to the character
$det(J/J^2) = det(g_-)^{-1}$ of $G$. $\eta$ generates  $A^G$ as a $U$-right module for  the universal enveloping algebra $U=U(Lie(A))$.

\bigskip\noindent
For the odd superderivations $D'_i$ dual to the $\tilde\theta_i \in g_-$ for $i=1,..,s$
define 
 $$ \kappa_A = D'_s \circ \cdots \circ D'_1(\eta)  \ \in \ A^G \ .$$ 
 Since $Lie(G)$ acts on 
$k\cdot D'_n \circ \cdots \circ D'_1$ by the character $det(g_-)$, it is easy to see that $\kappa_A$ is annihilated
by all right-invariant derivations $D'_X$ for $X$ in $Lie(G)$. 
Furthermore $\kappa_A = 1$ modulo $A^G\cap J= J^G$
or
$$ e_A^*(\kappa_A) = 1 \ .$$ 

%
\bigskip\noindent

\goodbreak
\bigskip\noindent
\centerline{\bf A global splitting}

\bigskip\noindent
The even derivation $D=D_\theta :A\to A$ defined by the Euler operator
$$ D(x) \ =\ \sum_{i=1}^s \ \theta_i \cdot D'_i(x) \ $$
obviously satisfies $D(A) \subset J$ (with notations as in the last section). Hence as a derivation $ D(J^\nu) \subset J^\nu$  for all $\nu \geq 1$.
The map 
$$ E_\nu : \quad J^\nu/J^{\nu +1} \to  J^\nu/J^{\nu +1} $$
 induced by $D$ is $B$-linear. So it suffices to compute $E_\nu$ on the basis elements $\tilde\theta_I$.
[For $x\in J^\nu$ and $a\in A$ use that $D(xa)=xD(a) +D(x)a=D(x)a$ mod $J^{\nu +1}$
and  $D(A)\in J$ implies $D(xa)=D(x)a \mod J^{\nu +1}$.] Therefore, as an immediate consequence of $ D(\theta_j) = \theta_j$ modulo
 $J^2\cap A^G$, this shows
$$  E_\nu(\theta_I) = \nu \cdot \theta_I \quad , \quad E_\nu = \nu \cdot id_{J^\nu/J^{\nu +1}} \  .$$
 
 \bigskip\noindent
 {\bf Lemma 9}. {\it For $char(k)=0$ or $char(k)> s$ the even derivation $D: A\to J$ induces an $k$-linear isomorphism
 $$ D: J\cong J \ .$$}
 {\it Proof}. For large enough $\nu$ we have $J^{\nu +1} =0$. The diagram
 $$ \xymatrix{ 0 \ar[r] & J^{\nu +1} \ar[r]\ar[d]^D & J^{\nu} \ar[r]\ar[d]^D & J^\nu/J^{\nu +1} \ar[r] \ar[d]^{\nu \cdot id} & 0 \cr    
 0 \ar[r] & J^{\nu +1} \ar[r] & J^{\nu} \ar[r] & J^\nu/J^{\nu +1} \ar[r] & 0 \cr    }$$
 commutes.  Hence by downward induction $D:J^\nu\to J^\nu$ is an $k$-linear isomorphism  for all $\nu \geq 1$ using the snake lemma.  QED

\bigskip
The kernel $$\tilde B= kernel(D: A \to A) \ $$ 
of the derivation $D$
is a $k$-subalgebra of $A$.
In the situation of the last lemma the snake lemma for
$$ \xymatrix{ 0 \ar[r] & J \ar[r]\ar[d]^D_\cong & A \ar[r]\ar[d]^D & B \ar[r] \ar[d]^0 & 0 \cr    
 0 \ar[r] & J \ar[r] & A \ar[r] & B \ar[r] & 0 \cr    }$$
implies that
the restriction of the quotient homomorphism $\pi : A \to B$ to $\tilde B\subset A$ is 
bijective. This inverse of the isomorphism $ \pi: \tilde B \ \cong \ B $  then defines a splitting of $\pi: A\to B$. Hence we get

\bigskip\noindent
{\bf  Splitting theorem}. {\it Suppose $char(k)=0$ or $char(k)>s$. Then $\pi: \ \tilde B\ \cong\ B$ is even and there exists 
 an isomorphism
of $k$-superalgebras}
$$  A \ \ = \ \ A^G   \otimes^\varepsilon  \tilde B \ \ \cong \ \  k[\theta_1,...,\theta_s] \otimes^\varepsilon \tilde B \ .$$  

\goodbreak
\bigskip\noindent
\centerline{\bf Supergroups}

\bigskip\noindent
An  affine algebraic group $G$ acts on its Lie algebra
$g_+$ by the adjoint representation.  Let $g_-$ be any finite
dimensional algebraic representation of $G$ over $k$ with action
denoted by $Ad_-$. Then $g_+$ acts
on $g_-$ by derivations $ad_-=Lie(Ad_-)$. Consider $G$-equivariant quadratic
maps $$ Q: g_- \to g_+ $$
 with respect to these  actions of $G$ 
(i.e. arising from a symmetric $k$-bilinear form on $g_-$ with
values in $g_+$). A triple ${\bf G} = (G,g_-,Q)$ as above will
be called a supergroup (over $k$) provided
$$  ad_-(Q(v))\ v=0 $$
holds for all $v$ in $g_-$.
An associated Lie algebra $Lie({\bf
G})$ considered as a $\Z_2$-graded Lie algebra structure 
is defined on $g_+\oplus g_-$ in the
obvious way by the Lie bracket induced by the group structure of
$G$, the action of $G$ on $g_-$ and the map $Q$ (super
commutator). See [DM], p.59.

\bigskip\noindent
{\bf Example 1}. If ${\bf G} =(G,g_-,Q)$ is a supergroup, then
also its connected component in the Zariski topology ${\bf G}^0 =
(G^0,g_-,Q)$.

\bigskip\noindent
{\bf Example 2}. Super-affine Hopf algebra $A$ define 
supergroups $$(Spec(A/J),g_-,Q)\ ,$$ where $Q$ is the
restriction of the Lie bracket on $g_-$ to the diagonal.

\bigskip\noindent
{\bf Example 3}. As a special of example 2 for a finite dimensional super vector space
$V=V_+\oplus V_-$ over $k$ the standard supergroup $Gl(V)$ is
defined by $G=Gl(V_+)\times Gl(V_-)$ together with $g_-=
Hom_k(V_+,V_-) \oplus Hom_k(V_-,V_+)$ and $Q(A\oplus B)= \{A,B\}$
for the super commutator $\{A,B\}=A\circ B+B\circ A$. Here we used
the obvious identification $Lie(Gl(V_\pm))=End_k(V_\pm)$.

\bigskip\noindent
{\it Center}. For a super group ${\bf G}=(G,g_-,Q)$ let the center $Z({\bf G})$
be the maximal central subgroup of $G$, which acts trivial on $g_-$. 

\bigskip\noindent
{\it Morphisms}. A homomorphism $(G,g_-,Q)\to (G',g_-',Q')$ between
supergroups is a pair $\Phi=(\phi,\varphi)$, where $\phi:G\to G'$
is a group homomorphism between algebraic groups over $k$ and
where $\varphi:g_-\to g_-'$ is a $k$-linear $\phi$-equivariant map
such that $Q'(\varphi(X))=Lie(\phi)(Q(X))$.

\bigskip\noindent
{\it Representations}. A representation $(V,\Phi)$ of a supergroup ${\bf G} =(G,g_-,Q)$
is a finite dimensional $k$ super vector space $V$ together with a
homomorphism of supergroups $\Phi: (G,g_-,Q) \to Gl(V)$. The
category of such representations, also denoted ${\bf
G}$-modules, is a $k$-linear abelian rigid (monoidal) tensor category
 $$ Rep_k({\bf G}) $$
with the forget functor $(V,\Phi)\mapsto V$ as a super fibre
functor. This fiber functor factorizes over the functor
$$ Lie: Rep_k({\bf G}) \to Rep_k(Lie({\bf G})) \ .$$
The category $Rep_k(Lie({\bf G}))$ of super representations of the
Lie superalgebra $Lie({\bf G})$ again is a $k$-linear abelian rigid (monoidal) tensor
category. Notation: Let $\sigma$ be an automorphism of the supergroup ${\bf
G}$. If $(V,\Phi)$ is a ${\bf G}$-module, then also
$(V,\Phi\circ\sigma)$.

\bigskip\noindent

\goodbreak 
\bigskip\noindent
\centerline{\bf An equivalence of representation categories }

\bigskip\noindent
Suppose $k=\C$.  Let ${\cal H}$
 be the opposite of the category of affine super Hopf algebras over $k$. Let ${\cal HC}$ be the category of supergroups ${\bf G}=(G,g_-,Q)$. Recall $G$ is an affine algebraic groups over $k$, and morphisms in ${\cal HC}$ are algebraic with respect to the first component of the triples. There is an obvious forget functor $$ {\cal H} \to {\cal HC}\ .$$
 
\bigskip\noindent 
There is a similar forget functor from the category  ${\cal H}_\infty$ of differentiable Lie supergroups (as in [DM]) to the category  ${\cal HC}_\infty$  of differentiable Harish Chandra triples. Objects now are ${\bf G}_\infty=(G_\infty,g_-,Q)$ for classical Lie groups $G_\infty$. According to [DM] p.79, [CF], [K] p. 232 this forget functor is a quasi-equivalence of categories in the $C^\infty$-case. Consider the following commutative diagram of forget functors 
$$ \xymatrix{ {\cal HC}   \ar[r] &   {\cal HC}_\infty \cr  
{\cal H}   \ar[r]\ar[u] &   {\cal H}_\infty \ar[u]_\sim \cr} $$
Since an algebraic morphism is 
determined by its associated $C^\infty$ map, the functor ${\cal H} \to {\cal HC}$ is faithful by going over the top of the diagram.
We now show  

\bigskip\noindent
{\bf Theorem 4}. {\it The functor  ${\cal H} \to {\cal HC}$ is fully faithful.}
 
 \bigskip
 This immediately implies theorem 3 or the equivalent

\bigskip\noindent
{\bf Corollary 6}. {\it For  a super-affine Hopf algebra $A$ over $k=\C$ with its associated supergroup $\bf G$  there exists a tensor-equivalence of algebraic tensor categories over $k$
$$ Rep_k(A) \sim Rep_k(\bf G)\ .$$}

\noindent{\it Proof of theorem}. 
For $A,A'$ in $\cal H$ with associated triples
$Y'=(Spec(B'),g'_-,Q')$ and  $Y=(Spec(B),g_-,Q) $ in ${\cal HC}$  and a morphism 
$$\Phi: Y' \to Y$$ in ${\cal HC}$ we have to
construct a homomorphism of super Hopf algebras $\Phi^*:A\to A'$ inducing $\Phi$. By the diagram above the corresponding differentiable morphism  $\Phi_\infty$ 
exists in ${\cal H}_\infty$. 

\bigskip
By construction  $\Phi_\infty$ is  \lq{reduced algebraic}\rq, i.e. the underlying
morphism of Lie groups $G'_\infty \to G_\infty$ is induced from an algebraic morphism  $\Phi_{red}:G' \to G$ between the underlying
reduced algebraic groups. Hence it suffices, if reduced algebraic morphisms $\Phi_\infty$ of
${\cal H}_\infty$ are induced from algebraic scheme morphisms $\Phi^*$,
The algebraic scheme morphism then automatically respects the additional structures comultiplication, antipode and augmention; this is obvious, since by assumption the  $C^\infty$ morphism $\Phi_\infty$ induced from it has this property.

\bigskip
To construct $\Phi^*$ from a reduced algebraic $\Phi_\infty$ consider its graph $\Psi_{\infty}=(id,\Phi_\infty)$ 
$$ (id,\Phi_\infty): (G_\infty,g_-,Q) \to  (H_\infty,h_-,Q_H) = (G_\infty,g_-,Q) \times (G'_\infty,g'_-,Q') \ ,$$
which again is reduced algebraic. By projection onto the second factor it suffices to show that $\Psi_{\infty}$ is algebraic. Thus it is enough to consider 
reduced algebraic morphisms $\Psi_\infty$ which are closed immersions. This means
that the underlying Hopf algebra morphism $$\xymatrix{\Psi_{red}^*: B \otimes^\varepsilon B' \ar@{->>}[r] & B \cr}$$ is surjective, and that
the map  
$Lie(H_\infty) \hookrightarrow Lie(G_\infty)$
induced by $\Psi_{\infty}$ is injective. 

\bigskip\noindent
{\it Construction of $\Phi^*$}. We may assume that $\Phi_\infty$ is a locally algebraic closed immersion. How to find $\Phi^*$? By
the  splitting theorem it suffices to find a right vertical ring homomorphism $\varphi:A^G \to (A')^G$ 
$$ \xymatrix@C-6mm{ A \ar@{.>}[d]^{\Phi^*} & \cong & \tilde B  \ar[d]^{\Phi_{red}^*} &  \otimes & A^G \ar@{.>}[d]^{\varphi} \cr
A' & \cong & \tilde B' & \otimes &  (A')^{G'} \cr} $$
such that the morphism of super schemes $\Phi^*$ induced on the left extends
to the given $\Phi_{\infty}$ in the differentiable category. Such $\varphi$ of course  exists if an only
if the pullback $\Phi^*_\infty$ of superfunctions in the differentiable sense satisfies the algebraicity condition
$$   \Phi^*_\infty(A^G)  \subset (A')^{G'} \ .$$   
Now use $Lie(H) = Lie(H_\infty)$ and $Lie(G) = Lie(G_\infty)$,
being defined by left-invariant derivations $D_X$ on the super ring of algebraic
resp. differentiable functions. For $X\in g'_+ \subset g_+$ there is 
a commutative diagram
$$ \xymatrix{   C^\infty(Y_{\infty})   \ar[r]^{D_X} \ar@{->}[d]_{\Phi_\infty^*} & C^\infty(Y_{\infty})     \ar@{->}[d]^{\Phi_\infty^*} \cr
C^\infty(Y'_{\infty})   \ar[r]^{D_X} &    C^\infty(Y'_{\infty}) \cr} $$
Since $ Lie(G'_\infty) \hookrightarrow Lie(G_\infty)$ the kernel $C^\infty(Y)^G$
of all $D_X$, $X\in g_+$ derivations on $C^\infty(Y)$ (being contained in the kernel
of all $D_X$ for $X\in g'_+$) pulls back
to the  kernel $C^\infty(Y')^{G'}$
of all $D_{X'}$, $X'\in g'_+$ on $C^\infty(Y')$.
Thus the desired existence of $\varphi$ is evident, if  the natural injection
$$   A'^{G'}  \hookrightarrow C^\infty(Y')^{G'} $$
is a bijection.  Notice $g'_+$ is a Lie algebra, hence integrable!  Thus $dim_k(C^\infty(Y')^{G'})= 2^s$ for $s=dim_k(g_-) $ as a consequence of the Frobenius theorem. See 
[DM], p.75   and [K], p. 230. 
Therefore $$dim_k(A'^{G'}) = 2^s= dim_k(C^\infty(Y')^{G'})\ .$$
This implies $A'^{G'} = C^\infty(Y')^{G'}$ and proves the claim. QED

\bigskip\noindent

\goodbreak
\bigskip\noindent
\centerline{\bf Semisimple tensor categories}

\bigskip\noindent
For a  $k$-linear abelian rigid (monoidal) tensor category $T$
with unit object $1_T$ and $End_T(1_T)=k$ the object $1_T$
is simple (see [DMi], prop 1.17). Furthermore

\bigskip\noindent
{\bf Lemma 10}. {\it  $T$ is semisimple iff
$1_T$ is injective or projective or $Hom_T(1_T,-)$ is exact or
$Ext^1_T(L,1_T)=0$ holds for all simple objects $L$ in 
$T$.}

\bigskip\noindent
{\it Proof}. $T$ is semisimple iff $Hom_T(N,M)=Hom_T(1_T,N^*\otimes M)=
Hom_T(N\otimes M^*,1_T)$ is exact
in $N,M$. This is equivalent to $Ext^1_T(L,1_T)=0$ for all (simple) objects $L$ in $T$. QED 

\bigskip\noindent
For tensor categories  $T$ and $T'$ as above
let $R: T \to T'$ be an exact covariant functor with an isomorphism
$\iota:  1_{T'} \cong R(1_T)$. Assume $ I: T' \to T$ is a
left-exact covariant functor. Let $p$ be  an epimorphism in $T$ 
$$ p: I(1_{T'}) \to 1_{T} \ .$$
Suppose there exists a natural transformation
$$ \nu: id \to R\circ I $$ 
such that $R(p) \circ \nu_{1_{T'}} = \iota$.

\bigskip\noindent
{\it Example}. $R$ exact  tensor functor with left adjoint $I$ . Then
$id\in Hom_{T}(I(W),I(W))$ defines $\nu_W\in Hom_{T'}(W,RI(W))$ 
and let $p\in Hom_{T}(I(1_{T'}),1_{T})$ correspond to $\iota \in Hom_{T'}(1_{T'},R(1_{T'}))$
for $\iota :1_{T'} \cong R(1_T)$. Then the above properties hold.

\bigskip\noindent
{\bf Lemma 11}. {\it a) In the situation above $T'$ is semisimple, if $T$ is semisimple.
b) If  $I$ is adjoint to $R$ and $End_{T'}(1_{T'})=k$, then $T$ is semisimple iff $T'$ is semisimple
and $p$ splits in $T$}.

\bigskip\noindent
{\it Proof}. b) Suppose $T'$ is semisimple. Then $Hom_{T}(I(1_{T'}),-) = Hom_{T'}(1_{T'},R(-))$ is exact. If  $p$ splits in $T$, then $  1_T \oplus I^+ = I(1_{T'}) $. Hence also $Hom_{T}(1_{T},-)$ is exact. Hence $T$ is semisimple. Conversely if $T$ is semisimple, $p$ splits.

\bigskip\noindent
a) Suppose $T$ is semisimple. If $T'$ is not semisimple, then
by the lemma 10 there exists a simple object $L$ and a nonsplit extension $E$ in $T'$
$$   0 \to 1_{T'} \overset{a}{\to} E \overset{b}{\to} L \to 0 \ .$$  
Since $\nu_{1_{T'}}: 1_{T'} \hookrightarrow RI(1_{T'})$ and $RI(a): RI(1_{T'}) \hookrightarrow RI(E)$
by our assumptions, $a({1_T'})\subset E$ is not in the kernel of $ \nu_E: E \to RI(E) $. Hence $b: kern(\nu_E)\to L$ is a monomorphism. 
Then  $kern(\nu_E)\neq 0$ implies $kern(\nu_E)\cong L$, since $L$ is simple. Since this would split $E$ this proves  $$\nu_E: E \hookrightarrow RI(E) \ .$$  Since $T$ is semisimple, $I(a): I(1_{T'} )\hookrightarrow I(E)$
 has a section $s: I(E) \to I(1_{T'})$. Then $$c: 1_{T'} \to R(1_T) $$ defined
by $  c= R(p) \circ R(s) \circ \nu_E \circ a  $ is nonzero.
[Otherwise $R(s)\circ R(I(a)) = id$, from $s\circ I(a)= id$,  
would give 
$ \iota =  R(p) \circ \nu_{1_{T'}} = R(p) \circ R(s) \circ R(I(a)) \circ \nu_{1_{T'}} = R(p) \circ R(s) \circ \nu_E \circ a = 0 $ by the naturality  $\nu_E \circ a = RI(a) \circ \nu_{1_{T'}}$  of $\nu$]. Hence  $c$ is an isomorphism as  $R(1_T)\cong 1_{T'}$ is simple, using $End_{T'}(1_{T'})=k$. Then $c^{-1}\circ R(p\circ s)$  splits $\nu_E(E)$
$$ \xymatrix{  0 \ar[r] & kernel \ar[r] & \nu_E(E) \ar[r]^{R(p\circ s)} & R(1_T) \ar[r] & 0 \cr
& & \nu_E(a(1_{T'})) \ar@{^{(}->}[u] \ar[ur]^\cong_c & & \cr }$$
 Since $\nu_E(E)\cong E$ this splits $a(1_{T'})$ in $E$. Contradiction! Hence $T'$ is semisimple.

\bigskip\noindent
 
\goodbreak
\bigskip\noindent
\centerline{\bf Semisimple representation categories}

\bigskip\noindent
 Let  ${\bf G}=(G,g_-,Q)$ be a supergroup over $k=\C$.  The obvious covariant exact restriction functor $ R : Rep_k({\bf G}) \to Rep_k(G) $ satisfies
$R(k)\cong k$. There exists a covariant induction functor $$ I : Rep_k(G) \to Rep_k({\bf G}) \ $$
which, for $V$ in $Rep_k(G)$ and $g= Lie({\bf G})$, is defined by $$ I(V) = U(g) \otimes^\varepsilon_{U(g_+)} V \ .$$  The action
of $g_+=Lie(G)$ on $I(V)$ comes from an algebraic action
of $G$ on $I(V)$ by the $g_+$-module isomorphism
$I(V)  \cong \Lambda^\bullet(g_-) \otimes^\varepsilon V$.
Hence $I(V)\in Rep_k({\bf G})$. It is easy to see that $I$ is exact and left adjoint to $R$, i.e. Frobenius reciprocity $Hom_{\bf G}(I(V),W) = Hom_G(V,R(W))$.  

\bigskip
Since $k$ has characteristic zero $Rep_k(G)$ is semisimple if and only if $G$ is a reductive algebraic group over $k$. Therefore lemma 11 b) implies

\bigskip\noindent
{\bf Theorem 5}. {\it  $Rep_k({\bf G})$ is semisimple 
if and only if (a) $G$ is reductive and (b) the surjection
of ${\bf G}$-modules defined by the adjunction morphism
$$ ad: I(k) \to k $$
has a splitting in the category $Rep_k({\bf G})$.}

\bigskip\noindent
{\it Remark}. By $char(k)=0$ condition a) holds iff
$g_+=Lie(G)$ is a reductive Lie algebra over $k$.
Condition b) says that the restriction $$   ad: I(k)^{\bf G} \to k $$ to the space of ${\bf G}$-invariant subspace of $I(k)$ is surjective.
For $g=Lie({\bf G})$ then $I(k)^{\bf G} = (I(k)^g)^G = (I(k)^g)^{\pi_0(G)}$ by [DG], prop. 2.1(c), p.309. 
The group of connected components $\pi_0(G)$ of $G$ in the Zariski topology is finite.
Since $char(k)=0$ the functor of $\pi_0(G)$-invariants is exact by Maschke's theorem. Hence conditions a) resp  b) are equivalent to the following conditions 
\begin{tabbing}
\quad \= \  a')\  \ The Lie algebra $g_+$ is reductive. \\
\quad \= \  b')\ \ The restriction of $ad: I(k)\to k$ to
$I(k)^g$ is surjective.
\end{tabbing}
{\bf Definition}. If $g$ satisfies these two properties, we say the  Lie superalgebra $g$  is {\it reductive}. $Rep_k({\bf G})$ is semisimple if and only if $g=Lie({\bf G})$ is reductive (by theorem 5). In this case we say ${\bf G}$ 
is reductive. 

\bigskip\noindent
{\it Connected component}.
As already explained \begin{tabbing} 
 \quad \= 1. \ ${\bf G}$ is reductive if and only if its connected component ${\bf G}^0$ is reductive.
\end{tabbing}
{\it Etale coverings}. Similarly we may replace $G$ by a finite etale covering
$G'\to G$. We say that the supergroup ${\bf G}'=(G',g_-,Q)$
attached to ${\bf G}=(G,g_-,Q)$ is a finite etale (central) cover of ${\bf G}$. Then of course
\begin{tabbing} \quad \= 2. \ ${\bf G}$ is reductive if and only if the etale cover ${\bf G}'$ is reductive.
\end{tabbing}
Then ${\bf G} = {\bf G'}/F$ for a finite subgroup $F$ of the center $Z({\bf G'})$ of ${\bf G}'$.

\bigskip\noindent
{\it Central quotients}. Finally if $Z$ is a closed subgroup of the center $Z({\bf G})$ of ${\bf G}$,
then ${\bf G}/Z = (G/Z,g_-,Q)$ again a supergroup called a central quotient.
Obviously 
 \begin{tabbing} \quad \= 3.
 \ If ${\bf G}$ is reductive, then any central quotient ${\bf G}/Z$ is also reductive
\end{tabbing}
since $Rep_k({\bf G}/Z)$ is a full subcategory of $Rep_k({\bf G})$, if $Rep_k({\bf G})$
is semisimple!

\bigskip\noindent

\goodbreak
\bigskip\noindent
\centerline{\bf Reductive supergroups}

\bigskip\noindent
The main classification statement
involves the orthosymplectic supergroups $$Spo(1,2r)=(Sp(2r,J),k^{2r},Q)\ .$$
Fix a nondegenerate antisymmetric $2r\times 2r$-matrix $J'=-J$ so that $g\in Sp(2r,k) \Leftrightarrow g'Jg=J$.
This identifies $sp(2r,J)$ with the matrices $X$ for which $JX$ is symmetric.
For the standard action $Ad_-$ of $Sp(2r,J)$ on $k^{2r}$ the  map $Q: k^{2r} \to sp(2r,J)$ $$ Q(v)_{\alpha\beta} = 
\sum_{\gamma=1}^{2r} v_\alpha v_\gamma J_{\gamma \beta} $$
for $v=(v_1,...,v_{2r})\in k^{2r}$ is  well defined and equivariant such that $Q(v)v=0$. So this defines a supergroup. Different choices of $J$ yield isomorphic supergroups.

\bigskip\noindent
{\bf Proposition 1}. {\it A supergroup is reductive over $k=\C$ if and only if its connected component
admits a finite etale central covering, which as a supergroup is a direct product of
super groups of the following type
\begin{tabbing}
\quad \=  1.  \ \  A classical central $k$-torus  \\
\quad \= 2 . \ \ Simple connected simply connected classical $k$-groups \\
\quad  \= 3. \ \ Simple supergroups of orthosymplectic type $Spo(1,2r)$ for integers $r\geq 1$.
\end{tabbing}}

\bigskip
Similarly a  Lie superalgebra is reductive if and only if, modulo the center, it is a direct sum
of simple  Lie superalgebras of classical type or of the orthosymplectic types $BC_r=spo(1,2r)$ corresponding
to the super groups $Spo(1,2r)$.

\bigskip\noindent
{\it Proof}. A product of reductive supergroups is reductive. We leave this as an exersize.
So in one direction it suffices that the supergroups $Spo(1,2r)$ are reductive.
In fact $Rep_k(Spo(1,2r))=Rep_k(spo(1,2r))$ because $Sp(2r,J)$ is simply connected, and this reduces to [DH], theorem 4.1. 

\bigskip
Now for the converse. By our preliminary remarks in the last section we may replace ${\bf G}^0$ by an etale finite covering ${\bf G}'$, where $G'=T\times S$ for  a $k$-torus  $T$ 
and a product $S$  of connected simple and simply connected $k$-groups. 
Then we can divide ${\bf G}'$ by its maximal central torus $Z$.
The new supergroup ${\bf G}''$ is reductive, if ${\bf G}$ is reductive. This
allows to reduce the proof to the case ${\bf G}=(G,g_-,Q)$ without central torus so that in addition $G$ is connected and a product of a torus $T$ and a simple simply connected $k$-group $S$.
If these conditions
hold and $Rep_k({\bf G})$ is semisimple, we say ${\bf G}$ is {\it good}.
So assume ${\bf G}$ is good. Then  by [DG], page 309ff and theorem 4 it suffices to prove that
$g=Lie({\bf G})$ is a product of simple  Lie superalgebras of the classical type
and types $BC_r$. Using condition b') this immediately would follow from [DH], theorem 4.1
for  semisimple $g_+$. 

\bigskip We already know $g_+$ is reductive.
To show that $g_+$ is semisimple
we claim that $g$ is a direct sum of  Lie superalgebras
$g_\nu$ 
with
 $(g_\nu)_+\neq 0$ and either  $(g_\nu)_-=0$  or  $(g_\nu)_-$ is an irreducible $(g_\nu)_+$-module with $(g_\nu)_+= [(g_\nu)_-,(g_\nu)_-]$.
This is easy: For $g_-=s\oplus t$ and an irreducible $g_+$-submodule $s$ $$ h_+=[s,s] $$ is an  ideal in $g_+$ by the Jacobi identity $[g_+,[s,s]] \subset [s,[g_+,s]] + [[g_+,s],s] \subset [s,s]$. Hence either $h=h_+$ in case $h_+$ commutes with $g_-$,
or otherwise
$$ h = h_+ \oplus s \ ,$$
is an ideal in $g$ with the desired property.
(As a ${\bf G}$-module, thus as a $g$-module)
  $g=h\oplus h'$ splits into
ideals by the semisimplicity of $Rep_k({\bf G})$. The ideal property $[h,h'] \subset h\cap h' =0$ decomposes  $g$. Since condition b') easily implies $h_+\neq 0$ for $s\neq 0$  (see [DH], prop.2.2) our claim follows by induction.

\bigskip
To show that $h=h_+ \oplus s$ is an ideal for $[h_+,g_-]\neq 0$,  notice $[h_+,t]=0$.
Indeed $[h_+,t]\subset [g_+,t] \subset t$ and $[h_+,t]=[[s,s],t]\subset [[s,t],s]\subset [g_+,s]\subset s$.
Thus $[h_+,s]=[h_+,g_-]\neq 0$. Therefore $s=[h_+,s] $, since $[h_+,s]$ is a $g_+$-submodule of $s$.
Obvious are  $[g_+,s]\subset s$ and $[g_+,h_+] \subset h_+$ and similarly
 $[g_-,h_+]=[g_-,[s,s]] \subset [[g_-,s],s] \subset [g_+,s] \subset s$.
To show $[g_-,s]=[t,s] + [s,s]\subset h_+$ 
use $[t,s]=[t,[h_+,s]]\subset [[h_+,t],s] + [[t,s],h_+] =
[[t,s],h_+] \subset [g_+,h_+] \subset h_+$.

\bigskip
If $(g_\nu)_- \neq 0$ is an irreducible
$(g_\nu)_+$-module, the center $z_\nu$ of $(g_\nu)_+$  acts by a character $\chi_\nu$. 
By the equivariance and surjectivity (!) of the Lie bracket $$ (g_\nu)_- \times (g_\nu)_- \to (g_\nu)_+ \neq 0 $$ the trivial action of $z_\nu$ on $(g_\nu)_+$ forces $2\chi_\nu=0$, hence $\chi_\nu=0$. Thus  $z_\nu$ is  in the center of $g$, therefore trivial
by our assumption that ${\bf G}$ is good. Hence the reductive Lie algebra $g_+$ is semisimple.
QED  

\bigskip\noindent

\goodbreak
\bigskip\noindent
\centerline{\bf The categories $Rep_k({\bf G},\varepsilon)$}

\bigskip\noindent

\bigskip\noindent
For a supergroup ${\bf G}=(G,g_-,Q)$ suppose $\varepsilon$ is in the center of $G(k)$
such that $\varepsilon^2=1$ and $Ad_-(\varepsilon)=-id_{g_-}$. Let
$ T=Rep_k({\bf G},\varepsilon)$ be the full subcategory of $Rep_k({\bf G})$
defined by the super representations $(V,\phi,\varphi)$ for which $\phi(\varepsilon)=\sigma_V$ is
the super parity automorphism $\sigma_V$ of $V$. $T$ is an algebraic tensor category
over $k$ (see [D]). 

\bigskip\noindent
Not every supergroup ${\bf G}=(G,g_-,Q)$ admits twisting elements $\varepsilon$
as above. But the extended supergroup ${\bf G}^{ext}=(G\times \mu_2,g_-,Q)$, where $Ad_-(g,\pm 1)=\pm Ad_-(g)$, always has the twisting 
element $\varepsilon^{ext} = (1,-1) \in G^{ext}=G\times \mu_2$. The forget
functor defines a tensor-equivalence
$$ Rep_k({\bf G}^{ext},\varepsilon^{ext}) = Rep_k({\bf G}) \ ,$$
since $(V,\phi,\varphi)\in Rep_k({\bf G})$ extends uniquely to $(V,\phi^{ext},\varphi)\in Rep_k({\bf G}^{ext},\varepsilon^{ext}) $
for $\phi^{ext}(g,\pm 1)= \sigma_V \phi(g) = \phi(g) \sigma_V$.

\bigskip\noindent
{\bf Lemma 12}. {\it $Rep_k(\bf G,\varepsilon)$ is semisimple if and only if
$Rep_k(\bf G)$ is semisimple.}

\bigskip\noindent
Since this is a statement on the underlying abelian categories, we may ignore the tensor structures on these categories. On the underlying abelian categories the parity change $\Pi(V) = V \otimes^\varepsilon \overline 1$, defined by the trivial super representation $\overline 1 = \Pi(1)$ on $k^{0\vert 1}$, induces
a functor
$ \Pi: Rep_k({\bf G}) \to  Rep_k({\bf G}) $
which in general does not preserve the subcategory $Rep_k({\bf G},\varepsilon)$.
However
$$ \Pi: Rep_k({\bf G}^{ext}) \to  Rep_k({\bf G}^{ext}) \ $$
preserves the subcategory $Rep_k({\bf G}^{ext},\varepsilon^{ext})$.

\bigskip\noindent
{\it Proof of lemma 12}. In the extended supergroup ${\bf G}^{ext}$
we have two twisting elements $\varepsilon$ and $\varepsilon^{ext}$. This
defines an element $z=\varepsilon\varepsilon^{ext} = (\varepsilon,-1) \in G\times \mu_2$ in the center of the supergroup ${\bf G}^{ext}$, i.e. $z$ is in the center of
$G^{ext}$ with trivial action on $g_-$,  and $z$ commutes with $\varepsilon$ and $\varepsilon^{ext}$. The eigenspace decomposition with respect to $z$ decomposes the category
$$ Rep_k({\bf G}^{ext}) =  Rep_k^+({\bf G}^{ext}) \oplus  Rep_k^-({\bf G}^{ext}) $$
and also its subcategories  $Rep_k({\bf G}^{ext},\varepsilon^{ext})$ and  $Rep_k({\bf G}^{ext},\varepsilon)$. Then by definition
$ Rep_k^+({\bf G}^{ext},\varepsilon^{ext}) = Rep_k^+({\bf G}^{ext},\varepsilon) $ and
$ Rep_k^-({\bf G}^{ext},\varepsilon^{ext}) = \Pi\bigl(Rep_k^+({\bf G}^{ext},\varepsilon^{ext})\bigr)$,
since $\varepsilon$ has trivial action and $\varepsilon^{ext}$ acts
by $-1$ on $\overline 1\in  Rep_k({\bf G}^{ext},\varepsilon^{ext})$.
Ignoring tensor structures
$ T= Rep_k({\bf G},\varepsilon) = Rep_k^+({\bf G}^{ext},\varepsilon) = Rep_k^+({\bf G}^{ext},\varepsilon^{ext}) $
and 
$ Rep_k({\bf G}^{ext},\varepsilon^{ext}) \ = \ Rep_k^+({\bf G}^{ext},\varepsilon^{ext}) \ \bigoplus\  \Pi(Rep_k^+({\bf G}^{ext},\varepsilon^{ext})) $ give
$$  Rep_k({\bf G}^{ext},\varepsilon^{ext}) =  T \bigoplus \Pi(T) \quad , \quad T=  Rep_k({\bf G},\varepsilon) \ .$$ 
Hence $T$ is semisimple iff  $Rep_k({\bf G}^{ext},\varepsilon^{ext}) =  Rep_k({\bf G})$
is semisimple. QED

\bigskip\noindent

\goodbreak
\bigskip\noindent
\centerline{\bf Remarks on ${\bf G}=Spo(1,2r)$}

\bigskip\noindent
We discuss the representations of the orthosymplectic group over $k=\C$.
The category $Rep_k({\bf G})$ of super representations of a supergroup ${\bf G}$
contains the  trivial even representation $1$ on $k=k^{1\vert 0}$
and the odd trivial representation $\overline 1$ on $k^{0\vert 1}$ such that
$\overline 1 \otimes^\varepsilon \overline 1 = 1$.  

\bigskip\noindent
For ${\bf G} =Spo(1,2r)$
the center of $G=Sp(2r,J)$ is $\mu_2$. The center of $\bf G$ is trivial. Hence $\varepsilon=-id$ gives a unique choice for a twisting element $\varepsilon$ to define a category
$T=Rep_k({\bf G},\varepsilon) \subset Rep_k({\bf G})$. Recall from the last section
$$ Rep_k({\bf G}) = T \ \bigoplus \ \Pi(T) \ .$$
Also notice $\Pi(W)=W\otimes^\varepsilon \overline 1$.

\bigskip\noindent
{\it The standard representation $V$}. Consider the  
following representation $(V,\phi,\varphi)\in Rep_k({\bf G},\varepsilon)$ of the supergroup
on ${\bf G}=Spo(1,2r)$. 
As a $G$-module $ V=V_+\oplus V_- = k \oplus g_-$ with trivial action on $V_+=k$ and with the standard representation of $G$ on $V_-$. This defines $\phi(X)\in End(V)_+$ for $X\in g_+$. We identify $V_-$ with $g_-$.   The odd elements  $v\in g_-$ act on $V$ by 
$\varphi(v)\in End(V)_-$ defined by the annihilation and creation operators
$$  \varphi(v) w = \frac{1}{2} v' J w \ \ \in V_+ \quad , \quad w \in V_- $$
$$  \varphi(v) \lambda = \lambda \cdot v \in V_-  \quad , \quad \lambda \in V_+\  .$$
Then $\phi(Q(v))= [\varphi(v),\varphi(v)]$ for $v\in g_-$. We call $V$ the orthosymplectic standard representation. It is easy to see that $V$ is an irreducible super representation.

\bigskip\noindent
{\it Invariant form $b$}. The orthosymplectic standard representation $V$ admits a nondegenerate supersymmetric ${\bf G}$-invariant form $$b: V\otimes^\varepsilon V \to k^{1\vert 0}$$
where $b$ is the orthogonal direct sum of the symmetric form $b(\lambda_1,\lambda_2) =
\lambda_1\lambda_2$ on $V_+=k$ and the antisymmetric form $b(v_1,v_2)=-\frac{1}{2}
v_1' J v_2$ on $V_-$. In fact the orthosymplectic supergroup ${\bf G}$ is the automorphism group of this supersymmetric form $b$ on $V$. In particular: The standard representation $V$
is an \lq{orthogonal self dual}\rq\ faithful representation of ${\bf G}$. Hence $V$ is a tensor generator of 
$$  T = Rep_k({\bf G},\varepsilon) = \ \langle V \rangle \ .$$
See [Sh] for an explicit decomposition of the tensor powers $V^{\otimes r}$.
See [RS] for a connection of $T$ with the representation category of the group $SO(2r+1)$.

\bigskip\noindent
{\bf Lemma 13}. {\it All irreducible representations in $T$ are \lq{orthogonal self dual}\rq.  All representations in $\Pi(T)$ are \lq{symplectic self dual}\rq.}

\bigskip\noindent
{\bf Proof}. If $W$ is \lq{orthogonal self dual}\rq\, then $\Pi(W)$ is \lq{symplectic self dual}\rq\ and vice versa. Since $Rep_k({\bf G}) = T \oplus \Pi(T)$ it therefore suffices that $T$ contains all \lq{orthogonal self dual}\rq irreducible representations. Tensor products of \lq{orthogonal self dual}\rq\  representations are \lq{orthogonal self dual}\rq, hence any multiplicity one representation contain in it is again \lq{orthogonal self dual}\rq.  By the  theory of highest weight vectors any irreducible representation $W$
in $Rep_k({\bf G})$ appears with multiplicity one in a tensor power of irreducible fundamental representations $V_i, i=1,..,r$ of ${\bf G}$ up to parity shift. For these
$(\overline 1^{\otimes i} \otimes^\varepsilon  V_i)_+ = \Lambda^i(g_-)$ and 
$(\overline 1^{\otimes i} \otimes^\varepsilon  V_i)_- = \Lambda^{i-1}(g_-)$.  See [Dj], p.31 and  p.36. Obviously $V_i \in Rep_k({\bf G},\varepsilon)$. The $V_i$ are self dual, therefore \lq{orthogonal self dual}\rq\ by considering
their restriction to $G$, which contains the highest weight representation with multiplicity
one as an \lq{orthogonal self dual}\rq\ representation of $G$. QED

\bigskip
We claim  

\bigskip\noindent
{\bf Lemma 14}. {\it For ${\bf G}=Spo(1,2r)$ the tensor subcategory of $Rep_k({\bf G})$ generated by the standard representation $V=k^{1\vert 2r}$
of ${\bf G}$  is $Rep_k({\bf G},\varepsilon)$. The tensor subcategory 
generated by $\Pi(V)$ is the full category $Rep_k({\bf G})$.}

\bigskip\noindent
{\it Proof}. It suffices to find $\overline 1 =\Pi(1) $ in a tensor power of $\Pi(V)$. Then $V=\Pi(V)\otimes^\varepsilon \overline 1$ generates $T$ and $T \ \bigoplus\ (T\otimes^\varepsilon \overline 1) = Rep_k({\bf G})$.
We claim $$\Pi( 1) \ \hookrightarrow \ \Pi \bigl(I(1)\bigr) \ \cong \ \Lambda^{2r+1}\bigl(\Pi(V)\bigr)\ $$
for the  induced module $I(k)=I(1)$. By Frobenius reciprocity  the dimension of $$End_{\bf G}\bigl(I(k)\bigr) \cong Hom_G\bigl(k,I(k)\bigr) \cong
\bigl(\Lambda^\bullet(g_-)\bigr)^G$$ is $r+1$ by the classical invariant theory of the group $G=Sp(2r)$. A basis for the invariants are the powers $\omega^i$ of the symplectic form $\omega\in \Lambda^2(g_-)$. Indeed
$$ I(k) = \bigoplus_{i=0}^r \ V_i \ \ \in \ \ Rep_k({\bf G},\varepsilon) $$
for  $ V_0 = 1 $ and  the different fundamental representations $ V_1, \cdots V_r $  
of ${\bf G}$ (see [Dj], p.36).
By Frobenius reciprocity also the dimension of
$$ Hom_{\bf G}\bigl(I(k),  \Lambda^{2r+1}\bigl(\Pi(V))\otimes^\varepsilon \overline 1 \bigr) = Hom_{G}\bigl(k, \Lambda^{2r+1}(\Pi(V))  \otimes^\varepsilon \overline 1\bigr)$$
is equal to $r+1$ using  $$\Lambda^{2r+1}\bigl(\Pi(V)\bigr)^G\ = \ \bigoplus_{j=0}^{2r+1} \ \Lambda^j(g_-)^G \ \otimes^\varepsilon \overline 1^{\otimes(2r+1-j) } \  \cong  \ \bigoplus_{j=0}^{r} \ \overline 1^{\otimes(2r+1-2j) } \ .$$ Then $I(k) \cong \Lambda^{2r+1}(\Pi(V))\otimes^\varepsilon \overline 1  $, provided $\Lambda^{2r+1}(\Pi(V))=\Pi( Sym^{2r+1}(V))$ has at least  
$r+1$ nonisomorphic irreducible constituents. 
For this (with the convention I of [DM], p.49 and p.62f) consider
the superpolynomial ring $S^\varepsilon(V)= Sym^\bullet(V)$ $$Sym^\bullet(V)= Sym^\bullet(V_+)  \otimes^\varepsilon \Lambda^\bullet(V_-) 
\ .$$  Multiplication with the invariant form $b\in Sym^2(V)$ is injective
inducing a filtration $  F_0 \hookrightarrow F_1 \hookrightarrow \cdots \hookrightarrow F_r$  of $F_r$ by ${\bf G}$-modules $ F_i = Sym^{2i+1}(V) \otimes^\varepsilon b^{\otimes (r-i)} $.
Notice $F_i \cong F_{i+1}$ for $i\geq r$. The highest weight submodules of the  $G_i=F_i/F_{i-1}$
for $i=0,..,r$ define $r+1$ nonisomorphic $G$-modules, 
since $(G_i)_- \cong
 \Lambda^{2i+1}(V_-)$ and $(G_i)_+ \cong
 \Lambda^{2i}(V_-)$ as $G$-modules. Hence  $$ I(1) \cong Sym^{2r+1}(V) \ .$$
Therefore the $G_i$ must have been irreducible ${\bf G}$-modules.
Considering highest weights a comparison shows $G_i\cong V_{2i+1}$ for $0\leq i < \frac{r}{2}$  and $G_{r-i} \cong V_{2i}$ for $0\leq i \leq \frac{r}{2}$. Hence all the representations $V_i$ for $i=0,..,r$ are constituents of the tensor power $V^{\otimes (2r+1)}$. QED

\bigskip\noindent
We remark that there is a dual filtration $F'_i = Sym^{2i}(V) \otimes^\varepsilon b^{\otimes (r-i)} $ 
on $Sym^{2r}(V)$ with $G'_i \cong G_{r-i}$, and again $I(1) \cong Sym^{2r}(V)$.

 \bigskip\noindent
{\bf Lemma 15}. {\it Let $T$ be a semisimple algebraic tensor category
over an algebraically closed field $k$ of characteristic zero. For a simple object $W\neq 0$ of $T$ the categorial rank $rk_k(W)$ does not vanish.}

\bigskip\noindent
{\it Proof}. Since $rk_k(W)=sdim_k(W)$, this follows from [Ka], p.619 formula (2.6) with $B(0,n)=spo(2n,1)$ in the notations of loc. cit. QED

\bigskip\noindent

\goodbreak
\bigskip\noindent
\centerline{\bf Structure Theorem}

\bigskip\noindent

\bigskip\noindent
Assume $k= \C$.
Then according to proposition 1 a connected reductive supergroup $\bf G$ is of the
form ${\bf G} = ({\bf G}' \times H)/F$ where ${\bf G}' = \prod_{r\geq 1} Spo(1,2r)^{n_r}$ 
is a product of orthosymplectic supergroups and where $H$ is a reductive algebraic $k$-group.
Since $F$ is a finite central subgroup of  ${\bf G}' \times H$ and since the center of ${\bf G}'$ is trivial, this implies $F\subset H$. Hence ${\bf G} = {\bf G}' \times H'$ for $H'=H/F$. Hence

\bigskip\noindent
{\bf Lemma 16}. {\it A connected reductive supergroup ${\bf G}$ is isomorphic to
a product ${\bf G}' \times H$ where $H$ is a reductive algebraic $k$-group and where
${\bf G}' =   \prod_{r\geq 1} Spo(1,2r)^{n_r}$ is a product of orthosymplectic supergroups.}

\bigskip\noindent
For ${\bf G}=Spo(1,2r)$ and $G=Sp(2r)$ one has $Aut(G)=G_{ad}$ and therefore
$Aut({\bf G})= G$. In other words, any automorphism of ${\bf G}$ is an inner automorphism
$Int(g)$ for a unique element $g\in G$. Let ${\bf G}$ be a reductive supergroup. Then
the group $\pi_0({\bf G}) = \pi_0(G)$ acts on ${\bf G}'$. For $\overline g \in \pi_0(G)$ 
we can choose a representative $g\in G$, by a suitable modification with an element in
$G' =   \prod_{r\geq 1} Sp(2r)^{n_r}$, such that $g$ acts by a strict permutation of the factors
on ${\bf G}'$. The group of such $g\in G$ defines a canonical subgroup $G_1\subset G$
such that $G_1 \cap G' = 1$. Hence $G_1\subset H$. Hence any $\overline g\in \pi_0(G)=\pi_0(H)$ has a representative in $G_1\subset H$. We get a canonical homomorphism
$$ p: H \to  \ \prod_{r\geq 1}\ \Sigma_{n_r} $$
into the product of symmetric permutation groups  $\Sigma_{n_r}$ whose kernel is $G_1$.
Conversely given such a homomorphism $p: H \to  \prod_{r\geq 1} \Sigma_{n_r} $
for a reductive algebraic $k$-group $H$ one can construct the semidirect product supergroup
${\bf G}  = {\bf G}' \vartriangleleft  H$ obtained from the permutation action of $H$ on ${\bf G}' =   \prod_{r\geq 1} Spo(1,2r)^{n_r}$. Obviously in our case therefore

\bigskip\noindent
{\bf Theorem 6}. {\it  Any reductive supergroup ${\bf G}$ over an algebraically closed field $k$ of characteristic zero is isomorphic to a semidirect product ${\bf G}' \vartriangleleft H$
of a reductive algebraic $k$-group $H$ with a product ${\bf G}' =   \prod_{r\geq 1} Spo(1,2r)^{n_r}$ of simple supergroups of $BC$-type, where the semidirect product is defined by an abstract  group homomorphism
$$ p: \pi_0(H) \to \ \prod_{r\geq 1}\ \Sigma_{n_r} \ .$$}

\bigskip\noindent

\bigskip\noindent
\goodbreak
\centerline{\bf References}

\bigskip\noindent
[CF] Cook J.-Fulp R., Infinite dimensional Super Lie groups,
arXiv: math-ph / 0610061 


\bigskip\noindent
[CP] Chari V.-Pressley A., A guide to Quantum groups, Cambridge university press (1998)

\bigskip\noindent
[D] Deligne P., Categories Tensorielles, Moscow Mathematical Journal, vol. 2, n. 2, p. 227 - 248 (2002)

\bigskip\noindent
[DMi] Deligne P.-Milne J.S., Tannakian categories, in Hodge cycles, Motives, and Shimura
Varieties, Springer Lecture Notes in Mathematics 900 (1982) 

\bigskip\noindent
[DM], Deligne P. -Morgan J.W., Notes on supersymmetry (following J.Bernstein), AMS (1999), in Quantum fields and strings, a course for mathematicians, vol. 1.

\bigskip\noindent
[DG], Demazure M.-Gabriel P., Introduction to Algebraic geometry and Algebraic groups, North-Holland (1980) 

\bigskip\noindent
[Dj] Djokovic D.Z., Symplectic 2-graded Lie Algebras, Journal of Pure and Applied Algebra 9,
25-38 (1976)

\bigskip\noindent
[DH] Djokovic D.Z. -Hochschild G., Semisimplicity of 2-graded Lie Algebras II, Illinois J. Math., 20 (1976)

\bigskip\noindent
[H] Hochschild G., Semisimplicity of 2-graded Lie Algebras, Illinois J. Math., 20 (1976) 

\bigskip\noindent
[Ka] Kac V.G., Representations of classical
Lie superalgebras, in Differential Geometric methods in Mathematical Physics II, Springer Lecture Notes in Mathematics 676 (1978)

\bigskip\noindent
[K] Kostant B., Graded Manifolds, graded Lie theory, and Prequantisation, in
Differential Geometric Methods in Mathematical Physics, Springer Lecture Notes 
in Mathematics 570, (1977)

\bigskip\noindent
[RS] Rittenberg V.-Scheunert M., A remarkable connection between the 
representations of the Lie superalgebra osp(1,2n) and the Lie algebra
o(2n+1), Comm. Math. Phys. 83 (1982),
no. 1, 1-9.

\bigskip\noindent
[SR] Saavedra Rivano, Categories Tannakiennes, Springer Lecture Notes in Mathematics 265, (1972)

\bigskip\noindent
[Sch] Scheunert M. , The Theory of Lie Superalgebras, Springer Lecture Notes in Mathematics
716, (1979) 

\bigskip\noindent
[Sh] C.L.Shader J.Korean Math Soc 36 (1999), no. 3, p. 593-607

\bigskip\noindent
[S] Sweedler M.S., Hopf algebras, Benjamin (1968)

\end{document}